\documentclass[a4pitaper,11pt]{amsart}
\usepackage{amsmath,amsthm,amssymb}
\usepackage[mathscr]{eucal}
 \usepackage{cite}
\usepackage{upgreek}
\usepackage[bookmarks,bookmarksnumbered]{hyperref}
\usepackage{enumerate}

\numberwithin{equation}{section}

\newcommand{\car}{\curvearrowright}

\theoremstyle{plain}
\newtheorem{main}{Theorem}
\newtheorem{mcor}[main]{Corollary}

\newtheorem{theorem}{Theorem}[section]

\newtheorem{lemma}[theorem]{Lemma}

\theoremstyle{definition}

\newtheorem{notation}[theorem]{Assumption}
\newtheorem{remark}[theorem]{Remark}

\begin{document}

\title[Cocycle superrigidity for profinite actions of irreducible lattices]
{Cocycle superrigidity for profinite actions \\ of irreducible lattices}

\author[D. Drimbe]{Daniel Drimbe}
\address{Department of Mathematics, KU Leuven, Celestijnenlaan 200b, B-3001 Leuven, Belgium}
\email{daniel.drimbe@kuleuven.be}
\thanks{D.D. was partially supported by NSF Career Grant DMS \#1253402 and by PIMS fellowship.}

\author[A. Ioana]{Adrian Ioana}
\address{Department of Mathematics, University of California San Diego, 9500 Gilman Drive, La Jolla, CA 92093, USA.}
\email{aioana@ucsd.edu}
\thanks{A.I. was partially supported by NSF Career Grant DMS \#1253402 and NSF FRG Grant \#1854074.}

\author[J. Peterson]{Jesse Peterson}
\address{Department of Mathematics, Vanderbilt University, 1326 Stevenson Center, Nashville, TN 37240, USA.}
\email{jesse.d.peterson@vanderbilt.edu}
\thanks{J.P. was supported in part by NSF Grant DMS \#1801125 and NSF FRG Grant \#1853989.}

\begin{abstract} 
Let $\Gamma$ be an irreducible lattice in a product of two locally compact groups and assume that $\Gamma$ is densely embedded in a profinite group $K$. We give  necessary conditions which imply that the left translation action $\Gamma\car K$ is ``virtually" cocycle superrigid: any cocycle $w:\Gamma\times K\rightarrow\Delta$ with values in a countable group $\Delta$ is cohomologous to a cocycle which factors through the map $\Gamma\times K\rightarrow\Gamma\times K_0$, for some finite quotient group $K_0$ of $K$. As a corollary, we deduce that any ergodic profinite action of  $\Gamma=\text{SL}_2(\mathbb Z[S^{-1}])$ is virtually cocycle superrigid and virtually W$^*$-superrigid, for  any finite nonempty set of primes $S$.

\end{abstract}

\maketitle

\section{Introduction and statement of main results}

The study of measure preserving actions of countable groups on standard probability spaces up to orbit equivalence  has witnessed an explosion of activity in the last 20 years (see the surveys \cite{Sh05, Po07b,Fu11, Ga10, Va10, Io13,Io18}). Recall that two probability measure preserving (p.m.p.) actions $\Gamma\car (X,\mu)$ and $\Delta\car (Y,\nu)$  are called {\it orbit equivalent} (OE) 
if there is an isomorphism of probability spaces $\theta:X\to Y$ such that $\theta(\Gamma\cdot x)=\Delta\cdot\theta(x)$, for almost every $x\in X.$  If, in addition, there is a group isomorphism $\delta:\Gamma\rightarrow\Delta$ such that $\theta(g\cdot x)=\delta(g)\cdot\theta(x)$, for every $g\in\Gamma$ and almost every $x\in X$, then the actions are called {\it conjugate}.
 
The theory of orbit equivalence was initiated by H. Dye, in connection with the theory of von Neumann algebras \cite{MvN36}.  He proved that any two ergodic p.m.p.\ actions of the group of integers $\mathbb Z$ are orbit equivalent \cite{Dy59}. In the early 1980s, this result was extended to show that all ergodic p.m.p.\ actions of infinite amenable groups are orbit equivalent \cite{OW80} (see also \cite{CFW81}). In contrast, it was shown in \cite{Ep07} that any non-amenable group has uncountably many pairwise non-orbit equivalent free ergodic p.m.p.\ actions (see also \cite{Hj05,GP05,Io11a} for results addressing various important classes of non-amenable groups).

Moreover, the non-amenable case revealed a striking rigidity phenomenon: within certain families of actions of non-amenable groups, orbit equivalence implies conjugacy. The first OE rigidity results were obtained by R. Zimmer for actions of higher rank lattices via his cocycle superrigidity theorem \cite{Zi84}. By building upon Zimmer's work, A. Furman proved the remarkable fact that generic ergodic p.m.p.\ actions $\Gamma\car (X,\mu)$ of higher rank lattices (e.g., SL$_n(\mathbb Z)\car\mathbb T^n$, for $n\ge 3$), are OE {\it superrigid} \cite{Fu99}: any free p.m.p.\ action  which is OE to $\Gamma\car (X,\mu)$ must be ``virtually" conjugate to it. Subsequently, numerous impressive OE superrigidity results were obtained in \cite{MS06,Po07a,Po08, Ki10,Io11b,PV11,Fu11,Ki11,PS12,Io17, TD20,CK15,Dr18a,GITD19,BTD18}. 

We highlight here the breakthrough work of S. Popa who used his deformation/rigidity theory to prove that Bernoulli actions 
of property (T)  and product groups  are OE superrigid \cite{Po07a,Po08}.
Popa derived  this result from his seminal cocycle superrigidity theorem asserting that  any cocycle for such an action with values in a countable (more generally, $\mathcal U_{\text{fin}}$, see Definition 2.5 in \cite{Po07a}) group is cohomologous to a group homomorphism \cite{Po07a,Po08}. Later on, techniques and ideas from deformation/rigidity theory were used to prove several additional cocycle superrigidity results in \cite{Io11b,PV11,Fu11,PS12,TD20,Io17,Dr18a, GITD19,BTD18}.
Notably, the second author proved a cocycle superrigidity theorem for ergodic profinite actions $\Gamma\car (X,\mu)$ of property (T) groups  \cite{Io11b}. This shows that any cocycle  for $\Gamma\car(X,\mu)$ with values in a countable group is virtually (i.e., after restricting to a finite index subgroup $\Gamma_0<\Gamma$ and an ergodic component  of $\Gamma_0$) cohomologous to a group homomorphism.  Soon after, this result was generalized  to compact actions in \cite{Fu11}.  More recently, completing an analogy with Bernoulli actions, 
D. Gaboriau, R. Tucker-Drob and the second author proved that separately ergodic profinite actions of product groups are  cocycle superrigid  in the above sense \cite{GITD19}. 

In this paper, we establish cocycle superrigidity results for profinite actions of a new class of groups that arise as irreducible lattices in products of locally compact groups. This is in part motivated by a question in \cite{Th11} asking whether profinite actions of groups with property $(\uptau)$, and in particular the irreducible lattices SL$_2(\mathbb Z[1/p])$ for prime $p$, are cocycle superrigid (see Remark \ref{Thomas}). Additional motivation is provided, in view of the analogy between existing results for Bernoulli and profinite actions, by a recent result in \cite{BTD18}, following earlier results in \cite{PS12}, showing that Bernoulli actions of lattices in products of locally compact groups are cocycle superrigid.

Before stating our main result, we need to recall some terminology. 
Let $G\car (X,\mu)$ be a p.m.p.\ action of a locally compact second countable (l.c.s.c.) group $G$ on a standard probability space $(X,\mu)$. A sequence $\{A_n\}_{n\in\mathbb N}$ of measurable subsets of $X$ is said to be {\it asymptotically invariant} if it satisfies $\underset{n\to\infty}{\text{lim}}\text{sup}_{g\in F} \mu(g\cdot A_n\triangle A_n)=0$, for every compact set $F
\subset G.$ The action $G\car (X,\mu)$ is  called {\it strongly ergodic} if any asymptotically invariant sequence $\{A_n\}_{n\in\mathbb N}$ is trivial, in the sense that $\underset{n\to\infty}{\text{lim}}\mu(A_n)(1-\mu(A_n))=0.$ 
For a l.c.s.c. group $H$, a measurable map $w:G\times X\rightarrow H$ is called a {\it cocycle} if for all $g,h\in G$, we have that $w(gh,x)=w(g,h\cdot x)w(h,x)$, for almost every $x\in X$. Two cocycles $w_1,w_2:G\times X\rightarrow H$ are {\it cohomologous} if there is a measurable map $\phi:X\rightarrow H$ such that for all $g\in G$, $w_1(g,x)=\phi(g\cdot x)w_2(g,x)\phi(x)^{-1}$, for almost every $x\in X$.
Finally, let $\Gamma$ be a lattice in $G$ and $m_{G/\Gamma}$ be the unique $G$-invariant Borel probability measure of $G/\Gamma.$ 
For a p.m.p.\ action $\Gamma\overset{\alpha}{\car} (Y,\nu)$, we denote by $\text{Ind}_{\Gamma}^{G}(\alpha)$, the associated {\it induced action} $G\car (G/\Gamma\times Y, m_{G/\Gamma}\times \nu)$ (see the beginning of Section \ref{section2} for the precise definition of induced actions). 


Our main result shows that, under certain 
strong ergodicity assumptions, any cocycle of a profinite action of a lattice in a product of locally compact 
groups  is virtually cohomologous to a homomorphism. 

\begin{main}\label{A}
 Let $\Gamma$ be a countable dense subgroup of a compact profinite group $K$ and consider the left translation action $\Gamma\overset{\alpha}{\car} (K,m_K)$, where $m_K$ denotes the Haar measure of $K$. 
Write  $K=\underset{\longleftarrow}{\emph{lim }}K_n$ as an inverse limit of finite groups $K_n$, and let $r_n: K\to K_n$ be the quotient homomorphism. 
Suppose that $\Gamma$ is a lattice in a product of two compactly generated l.c.s.c. groups $G=G_1\times G_2$. Assume that the restrictions of $\emph{Ind}_{\Gamma}^{G}(\alpha)$ to $ G_1$ and $G_2$ are strongly ergodic and ergodic, respectively.

Let $w:\Gamma\times K\to\Delta$ be a cocycle with values in a countable group $\Delta$.

Then there is an integer $n$ such that $w$ is cohomologous to a cocycle $w':\Gamma\times K\to \Delta$ of the form $w'=w_0\circ (\emph{id}_{\Gamma}\times r_n),$ for some cocycle $w_0:\Gamma\times K_n\to\Delta.$


\end{main}

Next, we discuss two consequences of Theorem \ref{A}.
By \cite[Theorem B]{Io11b} any ergodic profinite action of a finitely generated group $\Gamma$ that has property (T) (or, more generally, has an infinite normal subgroup with the relative property (T)) is virtually cocycle superrigid.
As a consequence of Theorem \ref{A}, we obtain the first class of residually finite groups $\Gamma$ not admitting infinite subgroups with the relative property (T) whose every ergodic profinite translation action is virtually cocycle superrigid.

\begin{mcor}\label{S-adic}
Let $\Gamma=\emph{SL}_2(\mathbb Z[S^{-1}])$, for a finite nonempty set of primes $S$. Let $\Gamma\car^{\alpha} (X,\mu)$ be an ergodic profinite p.m.p.\ action. Write $\alpha$ as an inverse limit of actions
$\Gamma\car^{\alpha_n} (X_n,\mu_n)$, with 
 $X_n$ finite, and denote by $r_n:X\rightarrow X_n$ the quotient map, for every $n$.
 
Let $w:\Gamma\times X\to\Delta$ be a cocycle with values in a countable group $\Delta$. 

 Then there is an integer $n$ such that $w$ is cohomologous to a cocycle $w':\Gamma\times X\to \Delta$ of the form $w'=w_0\circ (\emph{id}_{\Gamma}\times r_n),$ for some cocycle $w_0:\Gamma\times X_n\to\Delta.$ 
\end{mcor}

\begin{remark}\label{Thomas} 
In \cite[Theorem 1.1]{Th11}, S. Thomas proved that if $n\geq 2$ and $S,T$ are sets of primes, then the classification
problem for the $S$-local torsion-free abelian groups of rank $n$ is Borel reducible to
that for the $T$-local groups of rank $n$ if and only if $S\subset T$. The proof of this result relies on cocycle superrigidity results for profinite actions of $\Gamma=\text{SL}_n(\mathbb Z[1/p])$, for $p$ prime. If $n\geq 3$, the proof of \cite[Theorem 1.1]{Th11} uses \cite[Theorem B]{Io11b}, which can be applied as $\Gamma$ has property (T). On the other hand, if $n=2$, then $\Gamma$ does have property (T) and the proof of \cite[Theorem 1.1]{Th11} is much more complicated and relies on Zimmer's cocycle superrigidity theorem. 

In this context, S. Thomas asked if the cocycle superrigidity theorem of \cite{Io11b} holds for groups $\Gamma$ with property $(\uptau)$. Corollary \ref{S-adic} gives a partial positive answer to this question in the case when $\Gamma=\text{SL}_2(\mathbb Z[S^{-1}])$, for a finite set of primes $S$ (for which property $(\uptau)$ has been established in \cite{LZ89}). As explained in \cite[footnote on page 3700 and Remark B.3]{Th11}, the case when $S=\{p\}$ consists of one prime can be used to considerably simplify the proof of \cite[Theorem 1.1]{Th11} for $n=2$.
\end{remark}

By using standard arguments (see  \cite{Io11b}),  Theorem \ref{A} and Corollary \ref{S-adic} imply that the actions from their statements are virtually OE-superrigid. Our next result shows that the actions covered by Corollary \ref{S-adic} are moreover virtually W$^*$-superrigid.
Recall that a free ergodic p.m.p.\ action $\Gamma\car (X,\mu)$  is called  {\it W$^*$-superrigid} if any free ergodic p.m.p.\ action $\Delta\curvearrowright (Y,\nu)$ giving rise to an isomorphic von Neumann algebra, $L^{\infty}(X)\rtimes\Gamma\cong L^{\infty}(Y)\rtimes\Delta$, must be conjugate to it.
The first families of W$^*$-superrigid actions were discovered about 10 years ago \cite{Pe10,PV10,Io11c}. Since then, many other families of W$^*$-superrigid actions have been found (see the introduction of  \cite{Dr18b}).

By combining Corollary \ref{S-adic} with S. Popa and S. Vaes' work \cite{PV14} we obtain the following:

\begin{mcor}\label{C}
Let $\Gamma\car (X,\mu)$ be any action as in Corollary \ref{S-adic}. Let $\Delta\car (Y,\nu)$ be any free ergodic p.m.p.\ action of a countable group $\Delta$.

Then $L^{\infty}(X)\rtimes\Gamma\cong L^{\infty}(Y)\rtimes\Delta$ if and only if there exist finite index subgroups $\Gamma_0<\Gamma$ and $\Delta_0<\Delta$, a $\Gamma_0$-invariant measurable set $X_0\subset X$ and a $\Delta_0$-invariant measurable set $Y_0\subset Y$ such that 
\begin{itemize}
\item $\Gamma\car X$ is induced from $\Gamma_0\car X_0$,
\item $\Delta\car Y$ is induced from $\Delta_0\car Y_0$, 
\item $\Gamma_0\car X_0$ is conjugate to $\Delta_0\car Y_0$, and 
\item $[\Gamma:\Gamma_0]=[\Delta:\Delta_0]$.
\end{itemize}
\end{mcor}

Here, we say that an ergodic p.m.p.\ action $\Delta\car (Y,\nu)$ is  {\it induced from} an action $\Delta_0\car Y_0$ if $\Delta_0<\Delta$ is a finite index subgroup, $Y_0\subset Y$ is a $\Delta_0$-invariant measurable set and $\nu(gY_0\cap Y_0)=0$, for all $g\in\Delta\setminus\Delta_0$.


\section{Cocycle rigidity for induced actions}\label{section2}

The goal of this section is to prove cocycle rigidity results for induced actions of translation actions.

We start by recalling the construction of induced actions and cocycles. Let $\Gamma$ be a lattice in a l.c.s.c.\ group $G$.  Let $p:G/\Gamma\rightarrow G$ be a Borel map such that $p(g\Gamma)\in g\Gamma$, for every $g\in G$. Define a cocycle $c:G\times G/\Gamma\rightarrow \Gamma$ by letting $c(g,x\Gamma)=p(gx\Gamma)^{-1}gp(x\Gamma)$, for every $g\in G$ and $x\Gamma\in G/\Gamma$.

Let $\Gamma\car^{\alpha} (Y,\nu)$ be a p.m.p.\ action and $w:\Gamma\times Y\rightarrow\Delta$ be a cocycle, for a l.c.s.c. group $\Delta$.  Put $(\tilde Y,\tilde \nu):=(G/\Gamma\times Y,m_{G/\Gamma}\times\nu)$. Then the {\it induced action} $\tilde\alpha:=\text{Ind}_{\Gamma}^{ G}(\alpha)$ of $G$ on $(\tilde Y,\tilde\nu)$ is defined by the formula $\tilde\alpha(g)(x\Gamma,y)=(gx\Gamma,\alpha(c(g,x\Gamma))y)$, and the {\it induced cocycle} $\tilde w:G\times\tilde Y\rightarrow\Delta$ associated to $w$ is defined by the formula $\tilde w(g, (x\Gamma,y))=w(c(g,x\Gamma),y)$, for every $g\in G, x\Gamma\in G/\Gamma$ and $y\in Y$.

We also recall that if $G \car (X,\mu)$ is a p.m.p.\ action of a l.c.s.c. group $G$ on a probability space $(X,\mu)$ and $\Delta$ is a countable group, then the {\it uniform distance} between two cocycles $w_1,w_2:G\times X\to \Delta$ is given by 
$$
\text{d}(w_1,w_2)=\underset{g\in G}{\text{sup}}\;\mu(\{x\in X\mid w_1(g,x)\neq w_2(g,x)\}).
$$ 

The following result extends \cite[Lemma 2.1]{Io11b} to locally compact groups, with an identical proof.

\begin{lemma}[\!\!{\cite[Lemma 2.1]{Io11b}}] \label{L: cohom}
Let $G{\car} (X,\mu)$ be an ergodic p.m.p.\ action of a l.c.s.c. group $G.$ Let $\Delta$ be a countable group, and $w_1,w_2: G\times X\to \Delta$ be two cocycles such that $\emph{d}(w_1,w_2)<1/8$.

Then there is a measurable map $\phi:X\to\Delta$ such that for all $g\in G$, $w_1(g,x)=\phi(gx)w_2(g,x)\phi(x)^{-1}$ for almost every $x\in X.$ Moreover, $\mu(\{x\in X\mid\phi(x)=1_\Delta\})>3/4.$
\end{lemma}

\begin{notation}\label{N: general}  Throughout this and the next section, we assume the following setting:
\begin{itemize}
\item Assume that $\Gamma$ is a countable dense subgroup of a compact profinite group $K=\underset{\longleftarrow}{\text{lim }}K_n$, where $K_n$ is a finite group, for every $n$. Let $r_n:K\rightarrow K_n$ be the quotient homomorphism.
\item  Let $\Gamma\overset{\alpha}{\car} (K,m_K)$ be the left translation action. 
\item Assume that $\Gamma$ is a lattice in a product of two l.c.s.c. groups $ G=G_1\times G_2.$  
\item Let $G\overset{\tilde\alpha}{\car}(X,\mu)=(G/\Gamma\times K,m_{G/\Gamma}\times m_K)$ be the induced action $\tilde\alpha:=\text{Ind}_{\Gamma}^{ G}(\alpha).$ 
\item We define an action $K\overset{\sigma}{\car}(X,\mu)$ which commutes with $\tilde\alpha$ by letting $\sigma(t)(x\Gamma,y)=(x\Gamma,yt^{-1})$.
For simplicity, we will use the notation $zt^{-1}:=\sigma(t)z$, for every $z\in X$ and $t\in K$.
\end{itemize}
\end{notation}

The following theorem is the main result of this section.

\begin{theorem}\label{Th: untwist}
Let $ G_0$ be a closed subgroup of $ G $ such that the restriction of $\tilde\alpha$ to $G_0$ is ergodic. 
Let $w: G_0\times X \to\Delta$ be a cocycle  for the restriction of $\tilde\alpha$ to $G_0$ with values in a countable group $\Delta$. For every $t\in K$, define a cocycle $w_t: G_0\times X \to\Delta$ by letting $w_t(g,x)=w(g,x t^{-1})$.
Assume that $\emph{d}(w_t,w)<1/32$, for every $t$ in a neighborhood $V$ of the identity $1_K$ of $K$.

Then there is an integer $n$ such that $w$ is cohomologous to a cocycle $w': G_0\times X\to\Delta$ of the form $w'(g,x)=w_0(g,(\emph{id}
_{G/\Gamma}\times r_n)(x))$, for some cocycle $w_0:G_0\times (G/\Gamma\times K_n)\to\Delta.$

\end{theorem}

{\it Proof.} We follow closely the proof of \cite[Theorem 5.21]{Fu11}.
By Lemma \ref{L: cohom}, $w_t$ is cohomologous to $w$ for every $t\in V.$ Therefore, there is a measurable map $f_t:X\to\Delta$ such that for all $g\in G_0$, 
\begin{equation}\label{deform}
w_t(g, x)=f_t(gx) w(g, x)f_t(x)^{-1}, \text{ for almost every } x\in X.
\end{equation}
Moreover, the map $f_t$ satisfies $\mu (\{x\in X\mid f_t(x)=1_\Delta\})\geq 3/4$, for every $t\in V.$

Let $W\subset K$ be a neighborhood of $1_K$ such that we have $W^2\subset V$.
Let $t,s\in W$ and denote $F(x)=f_{ts}(x)^{-1}f_t(xs^{-1})f_s(x)$. Using \eqref{deform} twice, we obtain for all $g\in G_0$
$$F(gx)=w(g,x)F(x)w(g,x)^{-1}, \text{ for almost every }  x\in X.
$$ 
This implies that $F^{-1}(\{1_{\Delta}\})$ is $ G_0$-invariant and has positive measure since $f_s,f_t$ and $f_{ts}$ take the value $1_{\Delta}$ with probability at least $3/4$. Thus $\mu(F^{-1}(\{1_{\Delta}\}))=1$ by the ergodicity of ${\tilde\alpha}_{|G_0}.$ Therefore, for all $s,t\in W$, we have 
\begin{equation}\label{eq1}
f_{ts}(x)=f_t(x s^{-1})f_s(x), \text{ for almost every } x\in X.
\end{equation}


Since the finite index subgroups $L_n:=\ker(r_n)<K$ form a basis of neighborhoods of $1_K$, there is  $n$ such that $L_n\subset W$. 
We claim that there is a measurable map $\phi: X\to \Delta$ such that $$\text{$f_t(x)=\phi(x t^{-1})\phi(x)^{-1}$, for almost every $(t,x)\in L_n\times X.$}$$

To prove the claim, let $\mathcal R\subset K$ be a system of representatives for the left cosets of $K/L_n.$ Equation \eqref{eq1} shows that
 $f_{ts}(z,rx_0)=f_{t}(z,rx_0 s^{-1})f_{s} (z, rx_0),$
for all $r\in\mathcal R$, almost every $x_0,s,t\in L_n$ and almost every $z\in G/\Gamma$. By making the substitution $s=x_1^{-1}x_0$, we get that for all $r\in\mathcal R$, almost every $x_0,x_1,t\in L_n$ and almost every $z\in G/\Gamma$ we have
\begin{equation}\label{eq2}
f_t(z,rx_1)=f_{tx_1^{-1}x_0}(z,rx_0)f_{x_1^{-1}x_0}(z,rx_0)^{-1}.
\end{equation}
By Fubini's theorem we can find $x_0\in L_n$ such that \eqref{eq2} holds for all $r\in\mathcal R$, almost every $x_1,t\in L_n$ and almost every $z\in G/\Gamma$.
Define $\phi:X\to\Delta$ by letting $\phi(z,rx_1)=f_{x_1^{-1}x_0}(z,rx_0)$ for all $r\in\mathcal R,x_1\in L_n,z\in G/\Gamma$. Then \eqref{eq2} implies the claim.

Define the cocycle $w':G_0\times X\to \Delta$ by $w'(g, x)=\phi(gx)^{-1}w(g,x)\phi(x)$. Equation \eqref{deform} combined with the above claim implies that for every $g\in G_0$, we have $w'(g,xt^{-1})=w'(g,x),$ for almost every $t\in L_n$ and $x\in X$.  By Fubini's theorem we can find a map $w_0: G_0\times ( G/\Gamma\times K_n)\to \Delta$ such that for all $g\in G_0$, we have $w'(g,(x,y))=w_0(g,(x,r_n(y)))$, for almost every $(x,y)\in X.$ Then for all $g_1,g_2\in G_0$, $\tilde w_0$ satisfies the cocycle identity $w_0(g_1g_2,(x,y))=w_0(g_1,g_2 (x,y))w_0(g_2,(x,y))$, for almost every $(x,y)\in G/\Gamma\times K_n$. Moreover, if $q=\text{id}_{G_0}\times\text{id}_{G/\Gamma}\times r_n$, then $w_0\circ q=w'$, $(m_{G_0}\times m_{G/\Gamma}\times m_K)$-almost everywhere. Since $w'$ is measurable,  $w_0$ is measurable. Hence, $w_0$ is a cocycle, which finishes the proof. \hfill $\blacksquare$


The hypothesis of Theorem \ref{Th: untwist} requires that $w_t\rightarrow w$, as $t\rightarrow 1_K$, in the uniform metric. This assumption is guaranteed by the following lemma: 

\begin{lemma}\label{L: converge}
Assume that  the restriction of $\tilde\alpha$ to $G_1$  is strongly ergodic and $ G_2$ is compactly generated.
Let $w: G\times X\to \Delta$ be a cocycle into a countable group $\Delta$. For every $t\in K$, define a cocycle $w_t:  G\times X\to \Delta$ by letting $w_t(g,x)=w(g, xt^{-1})$.

Then $\emph{d}(w_{t| G_2},w_{| G_2})\to 0$, as $t\to 1_{K}.$

\end{lemma}
Here, $w_{| G_2}$ denotes the restriction of $w$ to $G_2\times X$.

The proof of this lemma follows closely \cite{GTD16} (see also the proof of \cite[Lemma 3.1]{GITD19}). Nevertheless, we include a detailed proof for the reader's convenience.

{\it Proof.} Let $S$ be a compact generating set for $ G_2$ and $\varepsilon\in (0,1)$.  Since the restriction of $\tilde\alpha$ to $G_1$   is strongly ergodic, there exist a compact set $F\subset  G_1$ and $\delta>0$ such that if $A\subset X$ is any measurable subset satisfying $\text{sup}_{g\in F}\mu (g^{-1}A\Delta A)<\delta$, then either $\mu (A)<\epsilon/4$ or $\mu (A)>1-\epsilon/4.$

For every $t\in K$ and $g\in G$, define $$A^t_g=\{x\in X\mid w(g,x)=w(g,xt^{-1})\}.$$
We claim that since $\Delta$ is countable, for every compact set $L\subset G$ we have that $\inf_{g\in L}\mu(A^t_g)\rightarrow 1$, as $t\rightarrow 1_K$.
To justify this, note that the formula $g\cdot (x,\rho)=(gx,w(g,x)\rho)$ defines a measure preserving  near action $G\curvearrowright (X\times \Delta,\mu\times c)$, where $c$ denotes the counting measure of $\Delta$.  Let $\pi:G\rightarrow\mathcal U(L^2(X\times\Delta))$ be the associated unitary representation. Since $\pi$ is measurable, it must be continuous (see, e.g., \cite[Theorem B.3]{Zi84}). Let $\tau:K\rightarrow\mathcal U(L^2(X\times\Delta))$ be the continuous unitary representation associated with the action $K\curvearrowright (X\times\Delta,\mu\times c)$ given by $t\cdot (x,\rho)=(xt^{-1},\rho)$.
Let $\xi={\bf 1}_{X\times 1_\Delta}\in L^2(X\times \Delta)$. Then $\|\xi\|_2=1$ and we have that $$\text{$ \pi(g)(\xi)={\bf 1}_{\{(x,\rho)\in X\times\Delta \mid w(g,g^{-1}x)=\rho\}}$, for every $g\in G$.}$$ This implies that $\mu(A_g^t)=\langle\tau(t)(\pi(g)\xi),\pi(g)\xi\rangle$, for every $g\in G$ and $t\in K$.
Since $\pi$ and $\tau$ are continuous in the strong operator topology, the set $\{\pi(g)\xi\mid g\in L\}\subset L^2(X\times \Delta)$ is compact and for every $g\in G$ we have that $\langle\tau(t)(\pi(g)\xi),\pi(g)\xi\rangle\rightarrow 0$, as $t\rightarrow 1_K$. This easily implies the claim.

 Next, the claim provides an open neighborhood $U\subset K$ of $1_K$ such for every $t\in U$ we have
\begin{equation}\label{almost}
\mu (A^t_g)>1-\delta/2 \text{ for all $g\in F$ and }\mu (A^t_h)>1-\epsilon/4 \text{ and for all $h\in S$}.
\end{equation}
Note that for any $t\in K, g\in  G_1$ and $h\in G_2$ we have
\begin{equation}\label{diff}
g^{-1}A^t_h\Delta A^t_h\subset X \setminus (A^t_g\cap h^{-1}A^t_g).
\end{equation}
Indeed, first notice that the cocycle relation implies
$
w(g,hx)w(h,x)= w(h,gx)w(g,x)$,  for almost every $x\in X.$ Now, if we take $x \in A^t_g\cap h^{-1}A^t_g$, then $w(g, x)=w(g, xt^{-1})$ and $w(g,hx)=w(g,hxt^{-1})$. Therefore, $x\in A^t_h$ if and only if $w(h,x)=w(h,xt^{-1})$ if and only if $w(h,gx)=w(h,gxt^{-1})$ if and only if $x\in g^{-1}A^t_h.$ This proves \eqref{diff}.

For every $g\in F, h\in G_2$ and $t\in U$, equations \eqref{almost} and \eqref{diff} imply that $\mu (g^{-1}A^t_h\Delta A^t_h)<\delta$. By our choice of $\delta$, it follows that 
\begin{equation}\label{erg}\text{ for all } t\in U \text{ and } h\in  G_2, 
\text{either } \mu (A^t_h)<\epsilon/4 \text{ or } \mu (A^t_h)>1-\epsilon/4.
\end{equation} 

We next claim that the set $ G_2':=\{h\in  G_2\mid \mu (A_h^t)>1-\epsilon/4, \text{ for all } t\in U \}$ is a subgroup of $ G_2.$ 
Note that for $A_{h_1}^t\cap h_1^{-1}A_{h_2}^t\subset A^t_{h_1h_2}$ for every $h_1,h_2\in  G_2$ and $t\in K.$ Therefore, if $h_1,h_2\in G_2'$ and $t\in U$, then $\mu (A^t_{h_1h_2})>\mu (A_{h_1}^t\cap h_1^{-1}A^t_{h_2})>1-\epsilon/2$. Since $1-\epsilon/2>\epsilon/4$, relation \eqref{erg} implies that $\mu (A^t_{h_1h_2})>1-\epsilon/4$, which implies $h_1h_2\in  G_2'.$  Since $A_{h^{-1}}^t=hA_h^t$ and thus $\mu(A_{h^{-1}}^t)=\mu(A_t^h)$, for every $g\in G$ and $t\in K$, this proves the claim.

Since $S\subset  G_2'$,  $ G_2'$ is a group and $S$ generates $G_2$, we get that $ G_2'= G_2.$ This proves the lemma.
\hfill $\blacksquare$

\section{Proof of Theorem \ref{A}}

We assume the setting from Assumption \ref{N: general}. Let $w:\Gamma\times K\to \Delta $ be a cocycle for $\alpha$ with values into a countable group $\Delta.$ 
Let  $\tilde w:G\times X\rightarrow\Delta$ be the induced cocycle for $\tilde\alpha$ defined by $\tilde w(g, (x\Gamma,y))=w(c(g,x\Gamma),y)$, for every $g\in G, x\Gamma\in G/\Gamma$ and $y\in K$.

For every $n$, recall that $r_n: K\to K_n$ denotes the quotient homomorphism and put $L_n:=\ker(r_n)$.
Define $X_n:= G/\Gamma\times K_n$ and $\tilde r_n:=\text{id}_{ G/\Gamma}\times r_n:X\to X_n.$

Since the restriction of $\tilde\alpha$ to $G_1$ is strongly ergodic,  Lemma \ref{L: converge} implies there is a neighborhood $V$ of $1_K$ in $K$ such that $\text{d}(\tilde w_{t|G_2},\tilde w_{|G_2})<1/32,$ for any $t\in V.$ 
Theorem \ref{Th: untwist} further implies that there exist an integer $n$, a map $\phi:X\to\Delta$ and a cocycle $v:  G_2\times X_n\to \Delta$ such that 
$$ \phi(gx)^{-1}\tilde w(g,x)\phi(x)=v(g,\tilde r_n(x)),
\text{ for all $g\in G_2$ and for almost every $x\in X.$}
$$

Define the cocycle $\tau: G\times X\to \Delta$ by
$\tau (g,x)=\phi^{-1}(gx)\tilde w(g,x)\phi (x)$ for $g\in  G$ and $x\in X.$
Note that
\begin{equation}\label{eq3}
\text{ $\tau (h,x)=v(h,\tilde r_n(x))$ for all $h\in G_2$ and a.e. $x\in X.$ }
\end{equation}
Therefore, for all $g\in G_1, h\in  G_2$ and for almost every $x\in X$, we have 
$$
\tau(g,hx)v(h,\tilde r_n(x))=v(h,\tilde r_n(gx))\tau(g,x),
$$
which is equivalent to
$$
\tau(g,hx)=v(h,\tilde r_n(gx))\tau(g,x)v(h,\tilde r_n(x))^{-1}.
$$
 For all $t\in K$ and $g\in G$ define $A^t_g=\{x\in X\mid \tau(g,x )=\tau (g,xt^{-1})\}.$ Since $L_n=\text{ker} (r_n)$, the set $A^t_g$ is $ G_2$-invariant, for every $g\in G_1,t\in L_n$. Using that $ \alpha_{| G_2}$ is ergodic, $A^t_g$ must be null or co-null, for every $g\in G_1,t\in L_n.$ 

Let $F\subset  G_1$ be a compact generating set. Then we can find  $N\geq n$ such that $A^t_g$ is non-null, for every $g\in F, t\in L_N$.
The previous paragraph then implies that $A^t_g$ is co-null, for every $g\in F,t\in L_N$. 
Since the set of $g\in  G_1$ such that  $A^t_g$ is co-null in $X$ for all $t\in L_N$ is clearly a subgroup of $G_1$ and $F$ generates $G_1$, we get that $A^t_g$ is co-null, for every $g\in G_1,t\in L_N$.
Moreover, \eqref{eq3} shows that $A^t_g$ is co-null, for every $g\in G_2,t\in L_N$. Since $ G= G_1\times  G_2$, it follows that $A^t_g$ is co-null, for every $g\in G,t\in L_N$.
This implies the existence of a cocycle $ \tau_0: G\times  X_N\to\Delta$ such that 
$$
\phi(gx)^{-1}\tilde w(g,x)\phi (x)=\tau(g,x)=\tau_0(g,\tilde r_N(x)), \text{ for all $g\in G$ and for almost every $x\in X$}.
$$

Finally, this together with Lemma \ref{L: induce}  below implies the existence of a cocycle $w_0:\Gamma\times K_N\to\Delta$ and a measurable map $\phi_0:K\to \Delta$ such that
$
w(g,x)=\phi_0(gx)w_0(g,r_N(x))\phi_0(x)^{-1}$ for all $g\in\Gamma$ and for almost every $x\in K.$ This finishes the proof. \hfill$\blacksquare$

We end this section with the following well-known result, whose proof we include for reader's convenience.

\begin{lemma}\label{L: induce}
Let $\Gamma$ be a lattice of a l.c.s.c. group $G$. Let $\Gamma\overset{\alpha}{\car} (X,\mu)$ and $\Gamma\overset{\beta}{\car} (X_0,\mu_0)$ be p.m.p.\ actions  such that there is a measurable onto map $\pi:X\to X_0$ satisfying $\pi(gx)=g\pi(x)$, for all $g\in \Gamma$ and almost every $x\in X.$
Denote by $G{\car} (\tilde X,\tilde\mu)$ and $G{\car} (\tilde X_0,\tilde\mu_0)$ the induced actions of $G$ associated to $\alpha$ and $\beta$, respectively.

Let $w:\Gamma\times X\to \Delta$ be a cocycle into a countable group $\Delta$ and denote by $\tilde w: G\times\tilde X\to \Delta$ the induced cocycle of $w$. Assume $\tilde w$ is cohomologous to a cocycle $\tilde w':G\times \tilde X\to\Delta$ of the form $\tilde w'(g,\tilde x)=\tilde w_0(g,(\emph{id}_{G/\Gamma}\times \pi)(\tilde x))$, for some cocycle $\tilde w_0: G\times \tilde X_0\to\Delta$.

Then $w$ is cohomologous to a cocycle $w':\Gamma\times X\to\Delta$ of the form $w'(\gamma,x)=w_0(\gamma,\pi(x))$, for some cocycle $w_0:\Gamma\times X_0\to\Delta$.
\end{lemma}

{\it Proof.} Let $p: G/\Gamma\to G$ be a Borel map such that $p(g\Gamma)\in g\Gamma$, for every $g\in G$. Define the cocycle $c:G\times G/\Gamma\to \Gamma$ by $c(g,z)=p(gz)^{-1}gp(z)$ for all $g\in G$ and every $z\in G/\Gamma$. Recall that for every $g\in G$ and $\tilde x=(z,x)\in \tilde X=G/\Gamma\times X$ we have $g\tilde x=(gz,c(g,z)x)$ and $\tilde w(g,\tilde x)=w(c(g,z),x).$

Let $\tilde \phi: \tilde X\to\Delta$ be a measurable map such that 
\begin{equation}\label{coh3}\tilde w(g,\tilde x)=\tilde\phi (g\tilde x) \tilde w_0(g,(\text{id}_{G/\Gamma}\times \pi)(\tilde x)) \tilde\phi (\tilde x)^{-1},\end{equation}
for all $g\in G$ and for almost every $\tilde x=(z,x)\in \tilde X$. 

By Fubini's theorem we can find  $z=h\Gamma\in G/\Gamma$, with $h\in G$, and a co-null subset $G_0\subset G$ such that for all $g\in G_0$, the identity \eqref{coh3} holds for almost every $x\in X$. Remark that there exists $g_0\in G$ such that $g_0 \Gamma h^{-1}\subset G_0.$ Indeed, this holds for any $g_0$ belonging to the co-null set $\cap_{\gamma\in\Gamma}G_0h\gamma^{-1}$.



For any $\gamma\in\Gamma$, denote $\gamma^h=g_0\gamma h^{-1}\in G_0$ and $\gamma_h=c(\gamma^h,h\Gamma)\in\Gamma$.
Note that $\gamma^h(h\Gamma,x)=(g_0\Gamma, \gamma_h x)$ and $\tilde w(\gamma^h,(h\Gamma,x))=w(\gamma_h,x)$, for every $x\in X$. Moreover, $\gamma_h=p(g_0\Gamma)^{-1}g_0\gamma h^{-1}p(h\Gamma)$.

Since $g^{-1}p(g\Gamma)\in\Gamma$ for all $g\in G$, it follows that the map $\gamma\to \gamma_h$ is a bijection of $\Gamma$. Hence, we can define a map $v_0:\Gamma\times X_0\to\Delta$ by letting $v_0(\gamma_h,x_0)=\tilde w_0(\gamma^h,(h\Gamma,x_0))$. Define also some measurable maps $\phi:X\to\Delta$ and $\psi:X\to\Delta$ by letting $\phi(x)=\tilde \phi (g_0\Gamma,x)$ and $\psi(x)=\tilde \phi (h\Gamma,x)^{-1}\tilde \phi (g_0\Gamma,x)$.
Therefore, for all $\gamma\in\Gamma$ and almost every $x\in X$ we have
$$
\begin{array}{rcl}
w(\gamma_h,x) &=& \tilde w (\gamma^h,(h\Gamma,x)) \\
&=&  \tilde\phi (\gamma^h(h\Gamma,x))\tilde w_0(\gamma^h,(\text{id}_{G/\Gamma}\times \pi)(h\Gamma,x))\tilde\phi (h\Gamma,x)^{-1}\\
&=& \tilde \phi ( g_0\Gamma,\gamma_h x)\tilde w_0(\gamma^h,(h\Gamma,\pi(x))) \psi (x) \tilde\phi(g_0\Gamma,x)^{-1} \\
&=& \phi(\gamma_h x) v_0(\gamma_h,\pi(x))\psi(x)\phi(x)^{-1}.
\end{array}
$$

Since the map $\gamma\to \gamma_h$ is a bijection of $\Gamma$, we obtain that $\phi(\gamma x)^{-1}w(\gamma,x)\phi(x)=v_0(\gamma,\pi(x))\psi(x)$ for all $\gamma\in\Gamma$ and almost every $x\in X$. By taking $\gamma=1_\Gamma$, if follows that $\psi(x)=v_0(1_\Gamma,\pi(x))^{-1}$. 
Define the map $w_0:\Gamma\times X_0\to \Delta$ by $w_0(\gamma,x_0)=v_0(\gamma,x_0)v_0(1_\Gamma,x_0)^{-1}$.
Therefore, the map $w':\Gamma\times X\to\Delta$ defined by $w'(\gamma,x)=w_0(\gamma,\pi(x))$ is a cocycle cohomologous to $w$. In particular, $w'$ is measurable, thus $w_0$ is measurable. Hence, $w_0$ is a cocycle, which finishes the proof.
\hfill$\blacksquare$

\section{Proof of Corollaries \ref{S-adic} and \ref{C}} Let $\Gamma=\text{SL}_2(\mathbb Z[S^{-1}])$, where $S$ is a finite nonempty set of primes.  Then $\Gamma$ is a lattice in $G=G_1\times G_2$, where $G_1=\text{SL}_2(\mathbb R)$ and $G_2=\prod_{p\in S}\text{SL}_2(\mathbb Q_p)$. 
Fix a positive integer $m$ with no prime factors from $S$. Denote $\Gamma(m):=\ker(\Gamma\rightarrow\text{SL}_2(\mathbb Z/m\mathbb Z))$ and consider the left translation action $G\car (G/\Gamma(m), \mu_m)$, where $\mu_m$ is the unique $G$-invariant Borel probability measure on $G/\Gamma(m)$. Let $\pi_m$ be the  associated Koopman unitary representation of $G$ on $L^2_0(G/\Gamma(m))=L^2(G/\Gamma(m))\ominus\mathbb C{\bf 1}_{G/\Gamma(m)}$.
Finally, let $$\pi=\oplus_{\{m\mid\; p\nmid m,\forall p\in S\}}\;\pi_m$$ be the direct sum of all such representations $\pi_m$.
The proof of Corollary \ref{S-adic} relies essentially on the following well-known fact.\begin{theorem}\label{gap}
The restrictions of $\pi$ to $G_1$ and $G_2$ have spectral gap.
\end{theorem}
 For the reader's convenience, we indicate below how this result follows from the literature.
Assuming Theorem \ref{gap}, we will now prove Corollary \ref{S-adic}.
\subsection*{Proof of Corollary \ref{S-adic}}
We will first prove the conclusion when $\Gamma\car^{\alpha} (X,\mu)$ is a left translation action $\Gamma\car \underset{\longleftarrow}{\text{lim }}\Gamma/\Gamma(m_n)$, for some sequence of positive integers $\{m_n\}$ containing no prime factors from $S$ and satisfying $m_n\mid m_{n+1}$, for all $n$. 
In this case, the induced action $\text{Ind}_{\Gamma}^G(\alpha)$ is isomorphic to the left translation action $G\car\underset{\longleftarrow}{\text{lim }}G/\Gamma(m_n)$. Since the Koopman representation of $G$ on $L^2_0(\underset{\longleftarrow}{\text{lim }}G/\Gamma(m_n))$ is isomorphic to a subrepresentation of $\pi$, Theorem \ref{gap}
implies that restrictions of $\text{Ind}_{\Gamma}^G(\alpha)$ to $G_1$ and $G_2$ are strongly ergodic.
Thus, in this case, the conclusion of Corollary \ref{S-adic} follows from Theorem \ref{A}.

In general, assume that  $\alpha$ is the inverse limit of a sequence of p.m.p.\ actions $\Gamma\car (X_n,\mu_n)$ with $X_n$ finite, for every $n$. Denote by $r_n:X\rightarrow X_n$ the $\Gamma$-equivariant quotient map.
Since $\alpha$ is ergodic, we may assume that $X_n=\Gamma/\Gamma_n$,  where $\{\Gamma_n\}_n$ is a descending chain of finite index subgroups of $\Gamma$. 

By a result of Serre \cite{Se70}, $\Gamma$ has the congruence subgroup property: any finite index subgroup of $\Gamma$ contains $\Gamma(m)$, for some positive integer $m$ having no prime factors from $S$.
Thus, we can find a sequence of positive integers $\{m_n\}$ such that $m_n$ contains no primes factors from $S$, $m_n\mid m_{n+1}$, and $\Gamma(m_n)\subset\Gamma_n$, for all $n$. 
Consider the profinite group $K=\underset{\longleftarrow}{\text{lim }}K_n$, where $K_n=\Gamma/\Gamma(m_n)$.   For every $n$, let $q_n:K\rightarrow K_n$ be the quotient homomorphism and denote $L_n=\ker(q_n)$.
Since the action $\Gamma\car (X,\mu)$ is a quotient of the left translation action  $\Gamma\car (K,m_K)$, we may identify it with the left translation $\Gamma\car (K/M,m_{K/M})$, for some closed subgroup $M<K$. 

Let $w:\Gamma\times K/M\rightarrow\Delta$ be a cocycle with values into a countable group $\Delta$. Define a cocycle $\tilde w:\Gamma\times K\rightarrow\Delta$ by letting $\tilde w(g,x)=w(g,xM)$.
Since the conclusion holds for the action $\Gamma\car K$ by the above,  $\tilde w$ is cohomologous to a cocycle which factors through $\text{id}_{\Gamma}\times q_{n_0}$, for some $n_0\geq 1$. Thus, we can find a homomorphism $\delta:\Gamma(m_{n_0})\rightarrow\Delta$ and a measurable map $\varphi:L_{n_0}\rightarrow\Delta$ such that \begin{equation}\label{coh}\text{$w(g,xM)=\varphi(gx)\delta(g)\varphi(x)^{-1}$, for all $g\in\Gamma(m_{n_0})$ and almost every $x\in L_{n_0}$. }\end{equation}

For $h\in L_{n_0}$, consider the set $S_h=\{x\in L_{n_0}\mid\varphi(xh)=\varphi(x)\}$.  Since $\lim\limits_{h\rightarrow 1_K}m_K(S_h)=m_K(L_{n_0})>0$, we can find $n_1\geq n_0$ such that $m_K(S_h)>0$, for all $h\in L_{n_1}$.
Now, if $h\in M\cap L_{n_0}$, then \eqref{coh} implies that $S_h\subset L_{n_0}$ is invariant under the left translation action of $\Gamma(m_{n_0})$.
Since $\Gamma(m_{n_0})<L_{n_0}$ is dense and $S_h$ is not null, we conclude that $m_K(L_{n_0}\setminus S_h)=0$, for all $h\in M\cap L_{n_1}$.
In particular, we have that $\varphi(xh)=\varphi(x)$, for almost every $x\in L_{n_1}$, for all $h\in M\cap L_{n_1}$. Thus, the restriction of $\varphi$ to $L_{n_1}$ factors through the quotient map $L_{n_1}\rightarrow L_{n_1}/(M\cap L_{n_1})$. 
Using the identification $L_{n_1}/(M\cap L_{n_1})\equiv L_{n_1}M/M$, it follows  that
  we can find a measurable map $\psi:L_{n_1}M/M\rightarrow\Delta$ such that $\varphi(x)=\psi(xM)$, for almost every $x\in L_{n_1}$.
By equation \eqref{coh} we thus have that 
\begin{equation}\label{coh2}
\text{$w(g,xM)=\psi(gxM)\delta(g)\psi(xM)^{-1}$, for all $g\in \Gamma(m_{n_1})$ and almost every $x\in L_{n_1}$.}
\end{equation}
This implies that $w$ is cohomologous to a cocycle $v:\Gamma\times K/M\rightarrow\Delta$ satisfying $v(g,x)=\delta(g)$, for all $g\in \Gamma(m_{n_1})$ and almost every $x\in L_{n_1}M/M$.  

We are now in position to apply an argument from the proof of \cite[Theorem B]{Io11b}. Let $\Gamma\car (X,\mu)$ be an ergodic profinite p.m.p.\ action, $\Gamma'<\Gamma$ a finite index subgroup, $X'\subset X$ a $\Gamma'$-ergodic component, $v:\Gamma\times X\rightarrow\Delta$ a cocycle and $\delta:\Gamma'\rightarrow\Delta$ a homomorphism such that $v(g,x)=\delta(g)$, for all $g\in\Gamma'$ and almost every $x\in X'$. Then parts 5 and 6 from the proof of \cite[Theorem B]{Io11b} show that there exists a finite $\Gamma$-invariant measurable partition $\{A_i\}_{i=1}^{\ell}$ of $X$ such that the map $v(g,\cdot):A_i\rightarrow\Delta$ is constant, for all $g\in\Gamma$ and $1\leq i\leq\ell$. 

Since $\Gamma(m_{n_1})<L_{n_1}$ is a dense subgroup, $L_{n_1}M/M\subset K/M$ is an ergodic component of the left translation action of $\Gamma(m_{n_1})$. 
The previous paragraph thus implies the existence of a finite $\Gamma$-invariant measurable partition $\{A_i\}_{i=1}^{\ell}$ of $X=K/M$ such that  the map $v(g,\cdot):A_i\rightarrow\Delta$ is constant, for all $g\in\Gamma$ and $1\leq i\leq\ell$.  
Moreover, \cite[Lemma 1.4]{Io11b} implies that we can find a positive integer $n$ such that $A_i$ is of the form $r_n^{-1}(Y)$, for some subset $Y\subset X_n$, for every $1\leq i\leq\ell$. This means that $v$ factors through the map $\text{id}_{\Gamma}\times r_n$, which finishes the proof.
\hfill$\blacksquare$


\subsection*{Proof of Theorem \ref{gap}} 
We will deduce this result from \cite{GMO08} by following closely the procedure from \cite[Section 6.3]{Lu94}.
Denote by $\mathcal P$ the set of all primes. 
 Let $H=\text{SL}_2(\mathbb R)\times(\prod'_{p\in\mathcal P}\text{SL}_2(\mathbb Q_p))$, where $\prod_{p\in\mathcal P}'\text{SL}_2(\mathbb Q_p)=\{(x_p)\in\prod_{p\in\mathcal P}\text{SL}_2(\mathbb Q_p)\mid \text{$x_p\in \text{SL}_2(\mathbb Z_p)$ for all but finitely many primes $p$}\}$ denotes the restricted product of $\text{SL}_2(\mathbb Q_p)$, $p\in\mathcal P$. Note that $H$ coincides with $\text{SL}_2(\mathbb A)$, where $\mathbb A$ is the Ad\`{e}le ring of $\mathbb Q$.
Consider the diagonal embedding of  $\Lambda=\text{SL}_2(\mathbb Q)$ into $H$. Then $\Lambda<H$ is a lattice. 
Consider the left translation action $H\car (H/\Lambda,\mu_\Lambda)$, where $\mu_\Lambda$ is the unique $H$-invariant Borel probability measure on $H/\Lambda$, and denote by $\rho$ the associated Koopman unitary representation of $H$ on $L^2_0(H/\Lambda)$.

Let $m$ be a positive integer with no prime factors from $S$. 
Write $m=p_1^{t_1}...p_k^{t_k}$, where $p_1,...,p_k\in\mathcal P\setminus S$ and $t_1,...,t_k\geq 1$. For every $p\in\mathcal P\setminus S$, we define an open subgroup $K_p<\text{SL}_2(\mathbb Z_p)$ as follows. If $p\notin\{p_1,...,p_k\}$, let $K_p=\text{SL}_2(\mathbb Z_p)$. If $p=p_i$, for $1\leq i\leq k$, let $K_p=\ker(\text{SL}_2(\mathbb Z_p)\rightarrow\text{SL}_2(\mathbb Z_p/p_i^{t_i}\mathbb Z_p))$. Then $K_m:=\prod_{p\in\mathcal P\setminus S}K_p$ is an open compact subgroup of $\prod'_{p\in\mathcal P\setminus S}\text{SL}_2(\mathbb Q_p)$. 

We claim that $GK_m\Lambda=H$.
By the Strong Approximation Theorem (see, e.g., \cite{LS03}) the diagonal embedding of $\Lambda$ into $\prod_{p\in\mathcal P}'\text{SL}_2(\mathbb Q_p)$ is dense.
This implies that $\text{SL}_2(\mathbb R)\Lambda$ is dense in $H$. Since $GK_m$ is an open subgroup of $H$ which contains $\text{SL}_2(\mathbb R)$, it follows that $GK_m\Lambda=(GK_m)(\text{SL}_2(\mathbb R)\Lambda)=H$. 

Since $G$ and $K_m$ commute, the subspace $L^2_0(H/\Lambda)^{K_m}\subset L^2_0(H/\Lambda)$ of $\rho(K_m)$-invariant vectors is $\rho(G)$-invariant.
On the other hand, $L^2_0(H/\Lambda)^{K_m}$ can be identified with $L^2_0({K_m}\backslash H/\Lambda)$
Since we have $H=GK_m\Lambda$, the $G$-space ${K_m}\backslash H/\Lambda$ is isomorphic to the $G$-space $G/(G\cap K_m\Lambda)$.  
Now, it is easy to see that $G\cap K_m\Lambda=\Gamma(m)$. By combining these facts, it follows that $\pi_m$, the Koopman representation of $G$ on $L^2_0(G/\Gamma(m))$, is isomorphic to the restriction of $\rho_{|G}$ to $L^2_0(H/\Lambda)^{K_m}$.

In conclusion, $\pi_m$ is isomorphic to a subrepresentation of $\rho_{|G}$, for every positive integer $m$ having no prime factors from $S$. Thus, $\pi$ is isomorphic to a subrepresentation of $\oplus_1^{\infty}\rho_{|G}$.
Now, by \cite[Theorem 1.11]{GMO08} there is some $s<\infty$ such that the positive definite function $H\ni g\mapsto \langle \rho(g)\xi,\eta\rangle\in\mathbb C$ belongs to $L^s(H)$, for all $\xi,\eta$ belonging to a dense subspace of $L^2_0(H/\Lambda)$. Hence, if $N\geq s/2$ is an integer, then $\rho^{\otimes_N}$ is contained is a multiple of the left regular representation of $H$. In combination with the above, we conclude that $\pi^{\otimes N}$ is contained is a multiple of the left regular representation of $G$. 
Hence, the restriction of $\pi$ to any non-amenable subgroup of $G$ has spectral gap. Since $G_1$ and $G_2$ are non-amenable, this implies the conclusion. \hfill$\blacksquare$

\subsection*{Proof of Corollary \ref{C}}  Let $\Gamma\curvearrowright (X,\mu)$ be a free ergodic profinite p.m.p. action of $\Gamma=\text{SL}_2(\mathbb Z[S^{-1}])$. As is well known, $\Gamma$ is measure equivalent to a direct product of $|S|$ non-abelian free groups (see, e.g., \cite[Remark 1.2]{DHI19} and the references therein). By \cite[Theorem 1.3]{PV14} we get that $\Gamma$ is Cartan rigid. Thus, if $\Delta\curvearrowright (Y,\nu)$ is any free ergodic p.m.p.\ action such that $L^{\infty}(X)\rtimes\Gamma\cong L^{\infty}(Y)\rtimes\Delta$, then the actions $\Gamma\curvearrowright (X,\mu)$ and $\Delta\curvearrowright (Y,\nu)$ are orbit equivalent. The conclusion of Corollary \ref{C} now follows from Corollary \ref{S-adic} through standard arguments (see the proofs of \cite[Theorem A]{Io11b} and \cite[Corollary C]{GITD19}). \hfill$\blacksquare$

\end{document}